\documentclass[11pt]{amsart}
\usepackage{amsthm,amsmath,amssymb}
\usepackage{mathtools}
\usepackage{graphicx}
\usepackage{caption}
\graphicspath{{Figure/}}
\numberwithin{equation}{section}
\numberwithin{figure}{section}
\theoremstyle{plain}
\newtheorem{definition}{Definition}[section]
\newtheorem{thm}{Theorem}[section]
\newtheorem{lem}[thm]{Lemma}

\newtheorem{prop}[thm]{Proposition}

\theoremstyle{definition}
\theoremstyle{remark}

\begin{document}
\setlength{\abovedisplayskip}{10pt}
\setlength{\belowdisplayskip}{10pt}

\title{A short proof of two shuffling theorems for tilings and a weighted generalization}

\author{Seok Hyun Byun}
\address{Department of Mathematics, Indiana University, Bloomington}
\email{byunse@indiana.edu}

\maketitle
\begin{abstract}
Recently, Lai and Rohatgi discovered a shuffling theorem for lozenge tilings of doubly-dented hexagons, which generalized the earlier work of Ciucu. Later, Lai proved an analogous theorem for centrally symmetric tilings, which generalized some other previous work of Ciucu. In this paper, we give a unified proof of these two shuffling theorems, which also covers the weighted case. Unlike the original proofs, our arguments do not use the graphical condensation method but instead rely on a well-known tiling enumeration formula due to Cohn, Larsen, and Propp. Fulmek independently found a similar proof of Lai and Rohatgi's original shuffling theorem. Our proof also gives a simple explanation for Ciucu's recent conjecture relating the total number and the number of centrally symmetric lozenge tilings.
\end{abstract}

\section{Introduction}
The enumeration of lozenge tilings of a region on a triangular lattice has received much attention during the last three decades. In particular, people tried to find regions whose number of lozenge tilings is expressed as a simple product formula. In [3], Ciucu considered an arbitrary string of triangles of alternating orientations that touch at corners and are lined up along a common axis. Then, he considered a hexagon with the triangles removed from its center and proved that the ratio of the number of lozenge tilings of two such regions is given by a simple product formula. In [4], Ciucu proved that the ratio of the number of centrally symmetric lozenge tilings of two such regions is also given by a simple product formula. In particular, he showed that the ratio for the centrally symmetric lozenge tilings is equal to the square root of the ratio for the total number of lozenge tilings. From this observation, he conjectured that this square root phenomenon holds in a more general setting.

Recently, Lai and Rohatgi [11] found and proved a shuffling theorem for lozenge tilings of doubly-dented hexagons, generalizing the work of Ciucu [3]. Later, Lai showed that similar theorems also exist for reflectively symmetric tilings [9] and centrally symmetric tilings [10]. In [10], he generalized the work of Ciucu and proved Ciucu's Conjecture in [4]. The proofs in [10] and [11] were based on Kuo's graphical condensation method (more precisely, Kuo's original recurrence from [8] for [11] and Ciucu's extension from [4] for [10]).

The goal of this paper is to give a unified and shorter proof for these two shuffling theorems, which also covers the weighted case. Unlike the original proofs, our arguments do not use the graphical condensation method but instead rely on a well-known tiling enumeration formula due to Cohn, Larsen, and Propp. Fulmek independently found a similar proof of Lai and Rohatgi's original shuffling theorem in [6]. Our proof also gives a simple explanation for Ciucu's recent conjecture relating the total number and the number of centrally symmetric lozenge tilings of regions with removed triangles.

\section{Statement of the theorem}

\begin{figure}
    \centering
    \includegraphics[width=8cm]{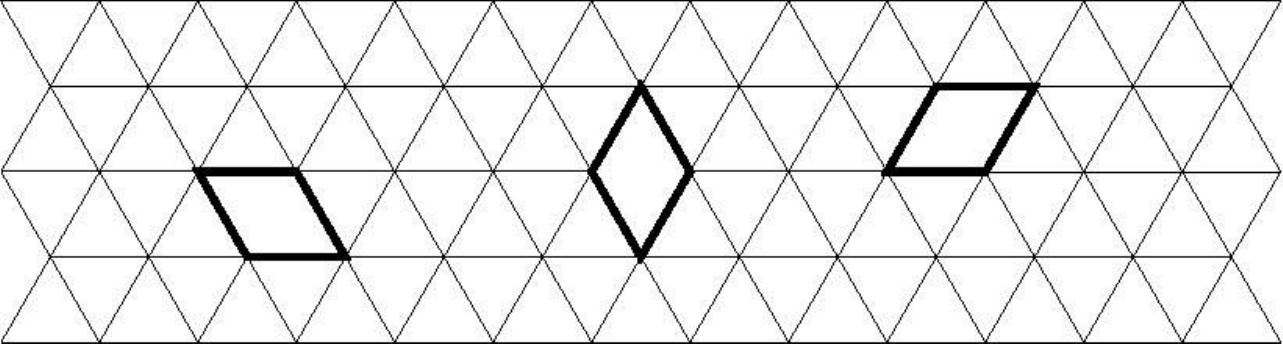}
    \caption{Left-, vertical-, and right-lozenge (from left to right).}
\end{figure}

In this paper, we are dealing with bounded regions on a triangular lattice and their lozenge tilings. A \textit{lozenge-shaped tile} (or \textit{lozenge} for short) is a union of two adjacent unit triangles on the triangular lattice, and a \textit{lozenge tiling} of a region is a collection of lozenges that covers the entire region without overlapping. We will always draw the lattice so that one family of lattice lines is horizontal. There are three kinds of lozenges that we can use. According to their orientation, we call them left-, vertical-, and right-lozenge, respectively (see Figure 2.1). For any region \(G\) on a triangular lattice, let \(M(G)\) be the number of its lozenge tilings. We now describe the region we are interested in. 

\begin{figure}
    \centering
    \includegraphics[width=10cm]{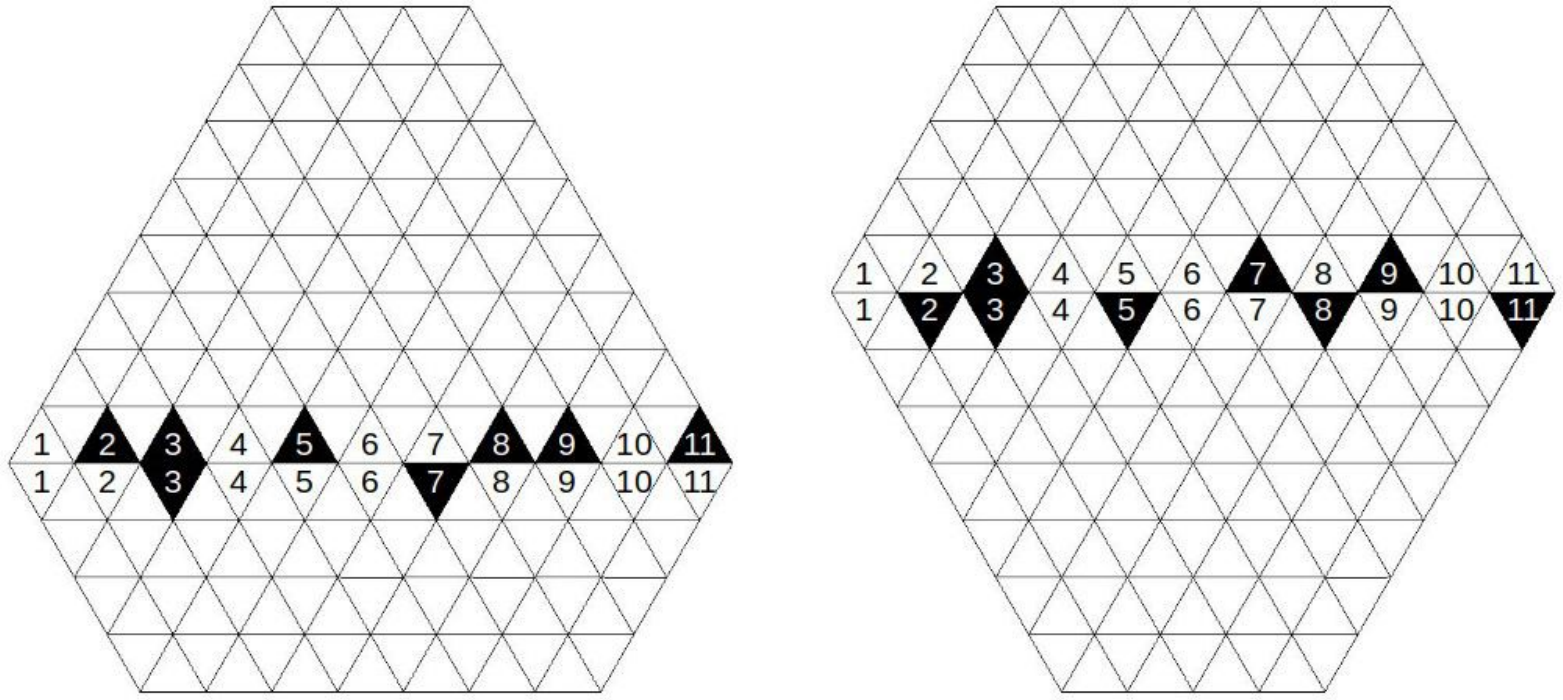}
    \caption{Two figures \(V_{3,8,4}(\{2,3,5,8,9,11\},\{3,7\})\) (left) and \(V_{6,5,7}(\{3,7,9\},\{2,3,5,8,11\})\) (right). The figure on the right is obtained from the left figure by flipping four of the removed up-pointing unit triangles (indexed by \(2,5,8\) and \(11\)), one of the removed down-pointing unit triangle (indexed by \(7\)) and moving up the horizontal diagonal \(3 (=4-1)\) units while preserving the height and width of the hexagon.}
\end{figure}

For nonnegative integers \(a,b\), and \(c\) such that \(c\leq a+b\), let \(V_{a,b,c}\) be the hexagon of side lengths \(a\), \(b\), \(c\), \(a+b-c\), \(c\), and \(b\) (clockwise from the top). Let \(\ell\) be the horizontal diagonal that connects the left vertex and the right vertex of the hexagon \(V_{a,b,c}\). This diagonal has length \(a+b\). Label the unit segments on it from left to right by \(1,\dotsc,a+b\).

For any subsets X and Y of \([a+b]\), where\([a+b]=\{1,\dotsc,a+b\}\), let \(V_{a,b,c}(X,Y)\) be the region obtained from \(V_{a,b,c}\) by removing up-pointing unit triangles whose bases are along the unit segments of \(\ell\) labeled by the elements of \(X\), and down-pointing unit triangles whose bases are along the unit segments of \(\ell\) labeled by the elements of \(Y\) (two examples are shown in Figure 2.2). 

\begin{definition}
For any finite subsets of the integers \(X\), \(X'\), \(Y\), and \(Y'\), we say that a pair \((X', Y')\) is a \textit{shuffling of} \((X, Y)\) if the following two conditions hold:

\begin{equation*}
\begin{aligned}
    &1) X\cup Y=X'\cup Y' \text{, and}\\
    &2) X\cap Y=X'\cap Y'.    
\end{aligned}
\end{equation*}
\end{definition}

In particular, if \(X\) and \(Y\) record the positions of the removed up- and down-pointing unit triangles along \(\ell\), and we are allowed to freely flip removed up-pointing unit triangles down and removed down-pointing unit triangles up, with the one restriction that pairs of removed unit triangles that form a vertical-lozenge are preserved, then the pair \((X',Y')\) recording the new positions of the removed unit triangles is a shuffling of \((X,Y)\).

In the paper [11], Lai and Rohatgi introduced weights on certain lozenges and the weighted enumeration of tilings under that weight. After introducing some notations, we also introduce similar weights on lozenges and consider the weighted enumeration of tilings\footnote{The weight that we now introduce is different from the one they used, but the argument of the current paper can be adapted to provide a short proof of their version.}.

For any set \(X=\{x_1,\dotsc,x_n\}\subset\mathbb{Z}_{+}\) of positive integers, where elements are written in increasing order, let \(\lambda(X)\) be the partition \((x_n-n,\dotsc,x_1-1)\), which may contain \(0\) as a part.
Also, for any finite disjoint subsets \(S, T\subset\mathbb{Z}_{+}\), let \(\Delta_{1}(S):=\displaystyle \prod_{s, s' \in S, s<s'}{(s'-s)}\), \(\Delta_{2}(S,T):=\displaystyle \prod_{s \in S ,t \in T}{|t-s|}\), \(\Delta_{1,q}(S):=\displaystyle \prod_{s, s' \in S, s<s'}{([s']_q-[s]_q)}\) and  \(\Delta_{2,q}(S,T):=\displaystyle \prod_{s \in S, t \in T, s<t}{([t]_q-[s]_q)}\cdot\prod_{s \in S ,t \in T, t<s}{([s]_q-[t]_q)}\), where \([n]_q:=\frac{1-q^n}{1-q}\) denotes the \(q\)-analogue of \(n\in\mathbb{Z}_{+}\). Also, let \(H(n):=\Delta_{1}([n])=\displaystyle\prod_{i=0}^{n-1}i!\).

\begin{figure}
    \centering
    \includegraphics[width=6cm]{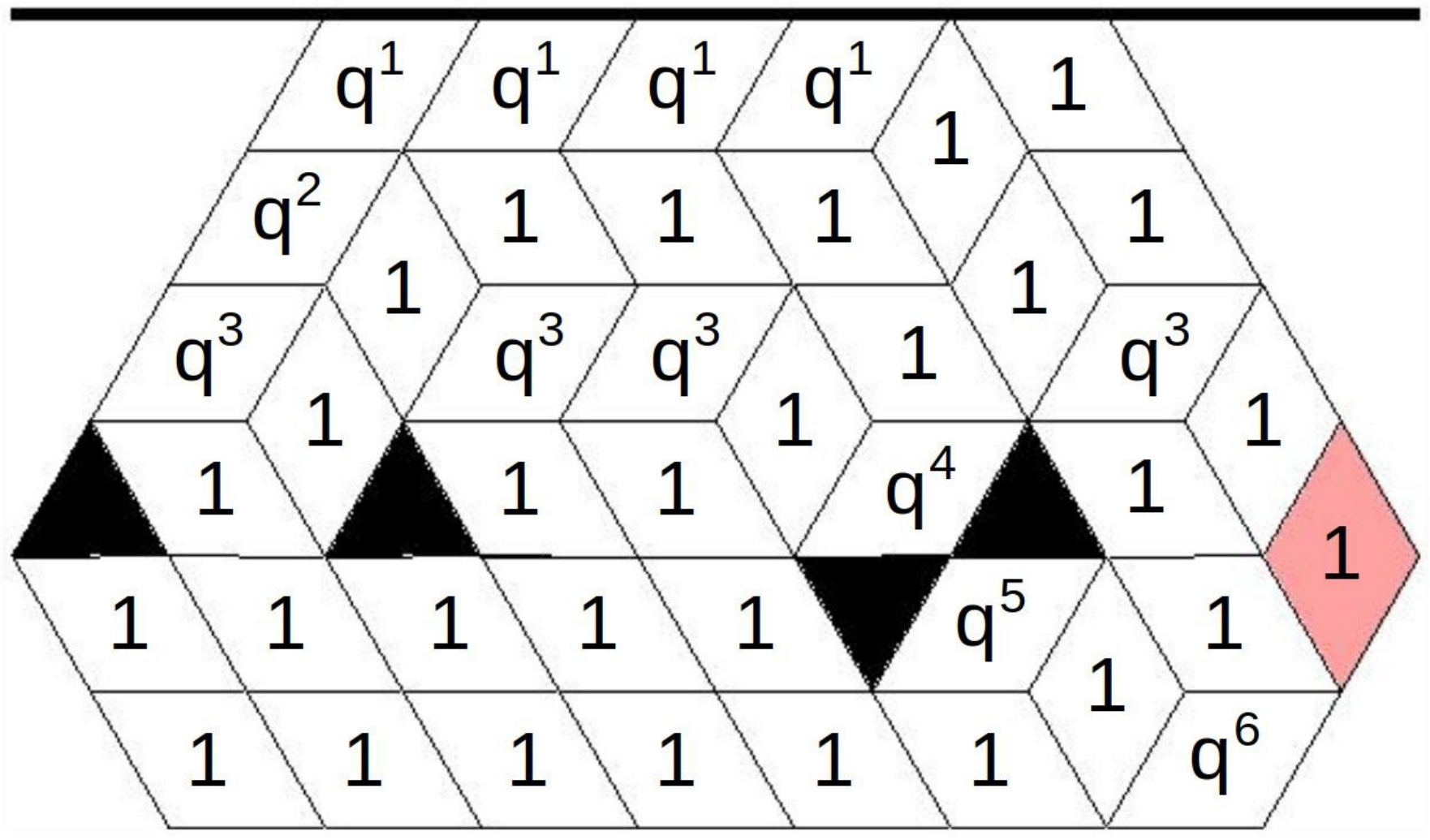}
    \caption{A tiling of \(V_{5,4,2}(\{1,3,7\},\{6\})\) with weighted lozenges. Its weight is \(q^{33}\).}
\end{figure}

\begin{figure}
    \centering
    \includegraphics[width=10cm]{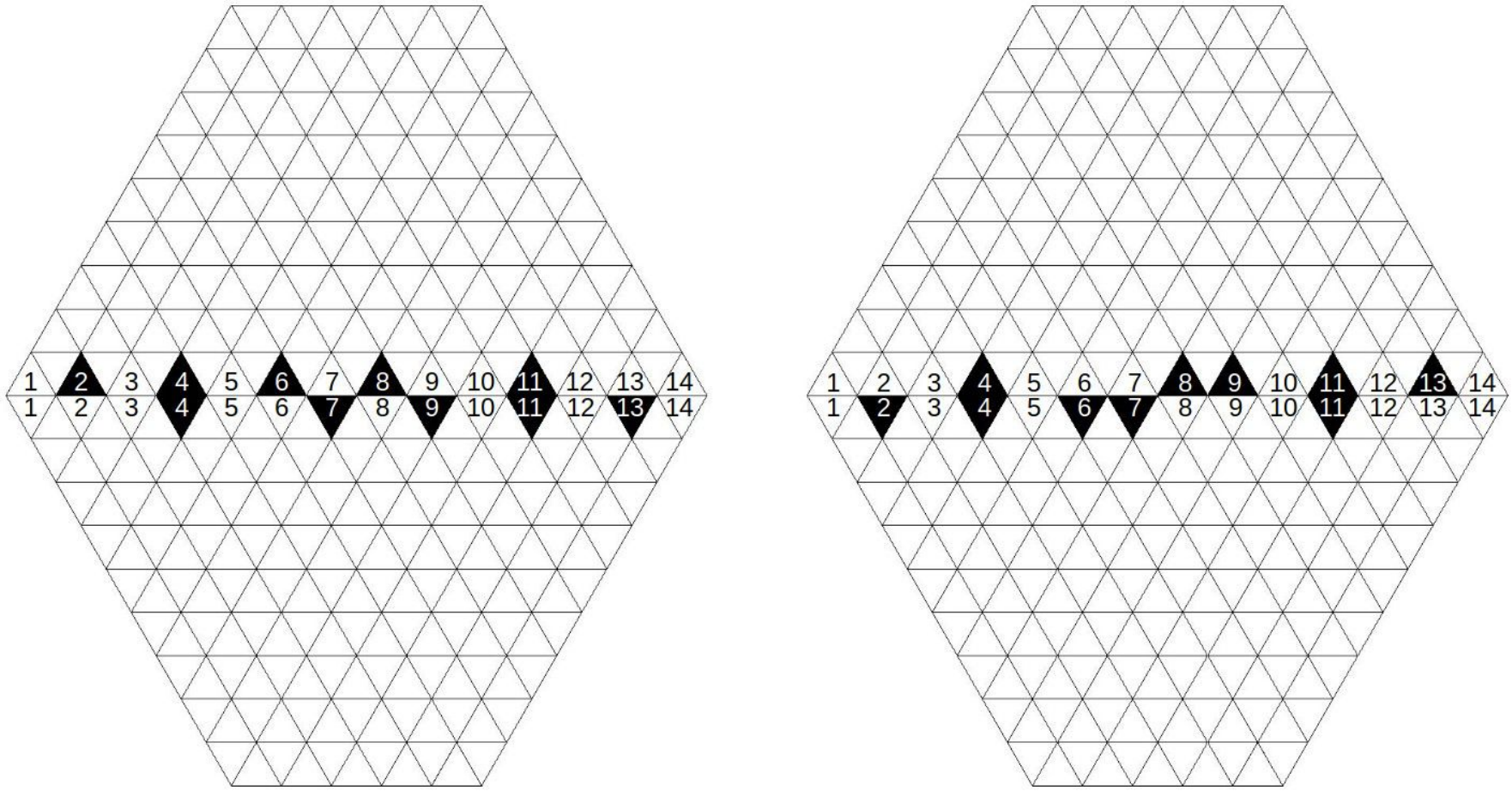}
    \caption{Two regions \(V_{5,9,9}(\{2,4,6,8,11\},\{4,7,9,11,13\})\) (left) and \(V_{5,9,9}(\{4,8,9,11,13\},\{2,4,6,7,11\})\) (right): Both regions are centrally symmetric. The region on the right is obtained from the left region by flipping two removed up-pointing unit triangles (indexed by \(2\) and \(6\)) and two removed down-pointing unit triangles (indexed by \(9\) and \(13\)).}
\end{figure}

For any bounded region \(G\) on the lattice, we give weight \(q^k\) to each right-lozenge (recall Figure 2.1) whose distance between the bottom side of the lozenge and the top side of \(G\) (= the highest horizontal line that intersects with the closure of the region \(G\)) is \(\frac{k\sqrt{3}}{2}\), and give weight 1 to all vertical- and left-lozenges (see Figure 2.3). When certain weights are given on lozenges and a lozenge tiling of the region is also given, the \textit{weight of the tiling} is the product of weights of all tiles that the given tiling contains. Also, the \textit{tiling generating function} of a region \(G\) is the sum of weights of tilings of \(G\), where the sum is taken over all lozenge tilings of the region \(G\). We will consider the tiling generating function under the weight described above and will denote it by \(M(G;q)\). Note that if we take \(q\rightarrow1\), then the tiling generating function \(M(G;q)\) becomes \(M(G)\).

A region \(G\) is \textit{centrally symmetric} if it is invariant under rotation by 180\(^{\circ}\). A lozenge tiling of a centrally symmetric region is \textit{centrally symmetric} if the tiling is invariant under rotation by 180\(^{\circ}\). For any centrally symmetric region \(G\), let \(M_\odot(G)\) be the number of its centrally symmetric lozenge tilings.
Also, for any positive integer \(k\), we say that two sets \(X\) and \(Y\subset[k]\) are \textit{k-symmetric} if \(Y=\{k+1-x|x\in X\}\) (or equivalently \(X=\{k+1-y|y\in Y\}\) ) holds, and denote this relation by \(Y=X_{(k)}\) (or \(X=Y_{(k)}\)). Obviously, if two sets are \(k\)-symmetric to each other, then they have the same cardinality. Note that the region \(V_{a,b,c}(X,Y)\) is centrally symmetric if and only if \(b=c\) and \(Y=X_{(a+b)}\) (two examples are shown in Figure 2.4).

Now we will state the two shuffling theorems. First part is a weighted generalization of the original shuffling theorem from [11], and the second part is the shuffling theorem for the centrally symmetric tilings from [10].

\begin{thm}
(a). Let \(a, b\), and \(c\) be nonnegative integers such that \(c\leq a+b\) and \(X,Y\) be subsets of \([a+b]\). Consider the region \(V_{a,b,c}(X,Y)\). While preserving removed unit triangles that form a vertical-lozenge, freely flip \(d\) of the removed up-pointing unit triangles along \(\ell\) down and \(u\) of the down-pointing up, so that their new positions are recorded by the sets \(X'\) and \(Y'\). Modify the boundary of \(V_{a,b,c}(X,Y)\) so that the height and width of the hexagon are preserved, but \(\ell\) is moved up \(d-u\) units (see Figure 2.2). This leads to the region \(V_{a',b',c'}(X',Y')\), where \((X',Y')\) is a shuffling of \((X,Y)\), and 

\medskip
\(\bullet\) \(a+b=a'+b'\) and \(b+c=b'+c'\)

\(\bullet\) \(a+b\geq c\) and \(a'+b'\geq c'\)

\(\bullet\) \(b-x=b'-x'\) and \(c-y=c'-y'\), where \(x, x', y,\) and \(y'\) are the cardinalities of \(X, X', Y,\) and \(Y'\), respectively.

\medskip
If $V_{a, b, c}(X,Y)$ has a lozenge tiling, then

\begin{equation}
    \frac{M(V_{a', b', c'}(X',Y');q)}{M(V_{a, b, c}(X,Y);q)}=q^{\alpha}\frac{\Delta_{1,q}([b])\Delta_{1,q}([c])}{\Delta_{1,q}([b'])\Delta_{1,q}([c'])}\frac{\Delta_{1,q}(X')\Delta_{1,q}(Y')}{\Delta_{1,q}(X)\Delta_{1,q}(Y)},
\end{equation}
where \(\alpha=\displaystyle \Bigg[\sum_{i'\in X'}i'-(b'+c')\sum_{j'\in Y'}j'-\frac{b'(b'+1)}{2}+(a'+b'+1)(b'+1)c'-\frac{(b'+1)c'(c'+1)}{2}+\frac{(a'+b'+1)c'(c'-1)}{2}\Bigg]-\Bigg[\sum_{i\in X}i-(b+c)\sum_{j\in Y}j-\frac{b(b+1)}{2}+(a+b+1)(b+1)c-\frac{(b+1)c(c+1)}{2}+\frac{(a+b+1)c(c-1)}{2}\Bigg]\).

(b). Let \(a\) and \(b\) be any nonnegative integers, and let \(X\) be a subset of \([a + b]\) of cardinality \(x\). Consider a centrally symmetric region \(V_{a, b, b}(X,X_{(a+b)})\). While preserving removed unit triangles that form a vertical-lozenge, freely flip removed up- and down-pointing unit triangles in \((a+b)\)-symmetric way (which means flipping removed up-pointing triangle whose position is indexed by \(i\) and removed down-pointing triangle whose position is indexed by \((a+b+1-i)\) at the same time), so that their positions are recorded by the sets \(X'\) and \(X'_{(a+b)}\) (see Figure 2.5). This leads to the region \(V_{a, b, b}(X',X'_{(a+b)})\), where \((X', X'_{(a+b)})\) is a shuffling of \((X, X_{(a+b)})\). If \(V_{a, b, b}(X,X_{(a+b)})\) has a centrally symmetric lozenge tiling, then
\begin{equation}
    \frac{M_\odot(V_{a, b, b}(X',X'_{(a+b)}))}{M_\odot(V_{a, b, b}(X,X_{(a+b)}))}= \sqrt{\frac{M(V_{a, b, b}(X',X'_{(a+b)}))}{M(V_{a, b, b}(X,X_{(a+b)}))}}=\frac{\Delta_{1}(X')}{\Delta_{1}(X)}.
\end{equation}

\end{thm}

\bigskip

We point out that given a region \(V_{a, b, c}(X,Y)\), the region \(V_{a', b', c'}(X',Y')\) is completely determined by the shuffling \((X',Y')\) of \((X,Y)\). Also, note that first equality in (2.2) implies Ciucu's conjecture on centrally symmetric lozenge tilings which says that the ratio between the numbers of centrally symmetric lozenge tilings of the two regions $V_{a, b, b}(X,X_{(a+b)})$ and $V_{a, b, b}(X',X'_{(a+b)})$ is exactly the square root of that of the plain counts of the lozenge tilings of the same regions. 

The rest of this paper is organized as follows. In section 3, we gather some previous results that we will employ in the proof of Theorem 2.1.
In section 4, we provide a proof of the theorem.

\section{Preparation for the proof}
The argument of the current paper was presented before in the author’s earlier work [2]. We present it here in a clearer form.

\medskip
From the assumption that the region \(V_{a,b,c}(X,Y)\) has a tiling, one can check that \(b-|X|=c-|Y|\) (It is because the region should contain the same number of up- and down-pointing unit triangles). Also, one readily sees that each tiling of \(V_{a,b,c}(X,Y)\) must contain precisely \(b-|X|\) (\(=c-|Y|\)) vertical-lozenges crossing the horizontal diagonal \(\ell\). Similarly, for \(V_{a',b',c'}(X',Y')\), the corresponding number is \(b'-|X'|\) (\(=c'-|Y'|\)). Note that the lengths of the horizontal diagonals in \(V_{a,b,c}(X,Y)\) and \(V_{a',b',c'}(X',Y')\) are the same. Furthermore,  
\(b-|X|=b'-|X'|\).

\medskip
Partition the set of tilings of \(V_{a,b,c}(X,Y)\) --- and also the set of tilings of \(V_{a',b',c'}(X',Y')\) --- in classes, according to the positions of the
\(b-|X|\) (\(=b'-|X'|\)) vertical-lozenges that straddle the diagonal \(\ell\). The proof will follow from the simple fact that the ratio of the cardinalities of corresponding classes in these two partitions is equal to a concrete simple product, which --- crucially --- turns out to be the \textit{same} for all classes of the partitions.

\medskip
This follows from the following result, which is the lozenge tilings interpretation given by Cohn, Larsen, and Propp [5] of a classical result due to Gelfand and Tsetlin [7]. 

\begin{figure}
    \centering
    \includegraphics[width=7.5cm]{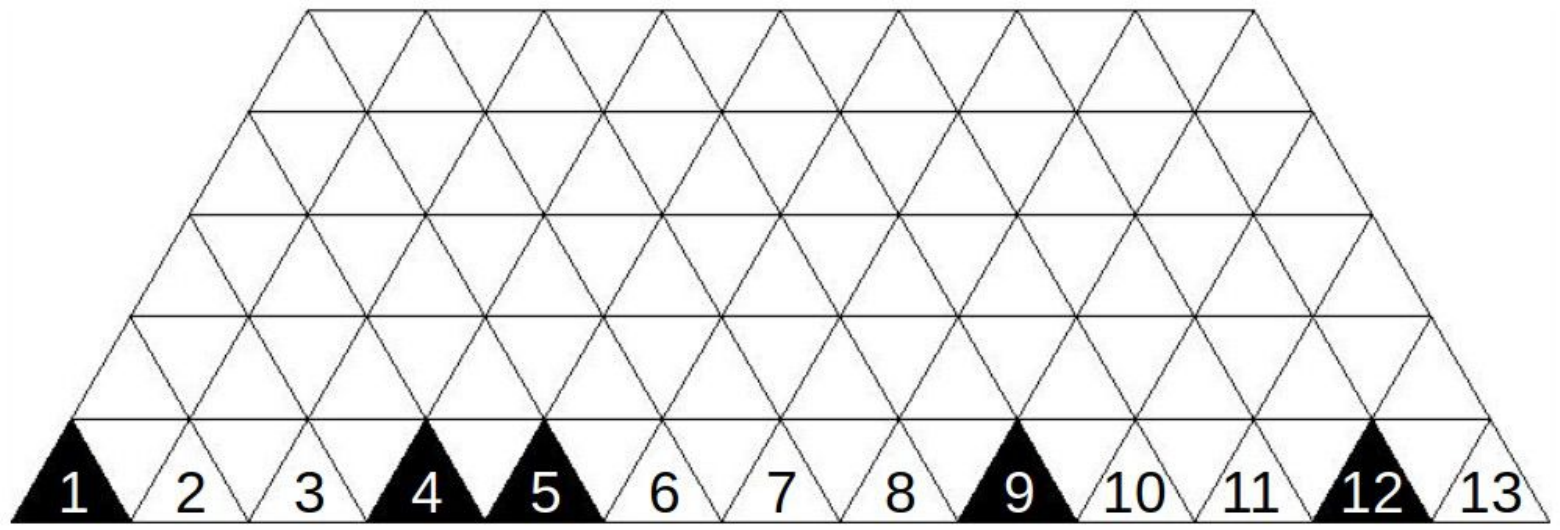}
    \caption{The region \(T_{8,5}(\{1,4,5,9,12\})\)}
\end{figure}

\begin{prop}
For nonnegative integers \(m\) and \(n\), let \(T_{m,n}\) be the trapezoid on the triangular lattice 
of side lengths \(m\), \(n\), \(m+n\), and \(n\) (clockwise from the top). Label the unit segments on the bottom from left to right by \(1,\dotsc,m+n\). 
For any subset \(S=\{s_1,\dotsc,s_n\} \subset [m+n]\), if \(T_{m,n}(S)\) is the region obtained from \(T_{m,n}\) by removing the up-pointing unit triangles whose bases have labels in \(S\) (see Figure 3.1 for an example), then
\begin{equation}
    M(T_{m,n}(S))=\frac{\Delta_{1}(S)}{\Delta_{1}([n])}=\frac{\Delta_{1}(S)}{H(n)}.
\end{equation}
\end{prop}

The same idea can be adapted to provide the proof of the weighted case. To deal with the tiling generating functions of the regions appeared in the theorem, we need the following Schur function identity.

\begin{lem}
Let \(X\), \(X'\), \(Y\), and \(Y'\) be any sets of positive integers whose cardinalities are \(x, x', y\), and \(y'\), respectively, so that a pair \((X', Y')\) is a shuffling of \((X, Y)\). If \(Z\) is any finite set of positive integers disjoint from \(X\cup Y(=X'\cup Y')\) whose cardinality is \(z\), then we have
\begin{equation}
\begin{aligned}
    &\frac{s_{\lambda(X'\cup Z)}(1,q,\dotsc,q^{x'+z-1})s_{\lambda(Y'\cup Z)}(1,q,\dotsc,q^{y'+z-1})}{s_{\lambda(X\cup Z)}(1,q,\dotsc,q^{x+z-1})s_{\lambda(Y\cup Z)}(1,q,\dotsc,q^{y+z-1})}\\
    &=\frac{\Delta_{1,q}([x+z])\Delta_{1,q}([y+z])}{\Delta_{1,q}([x'+z])\Delta_{1,q}([y'+z])}\frac{\Delta_{1,q}(X')\Delta_{1,q}(Y')}{\Delta_{1,q}(X)\Delta_{1,q}(Y)},
\end{aligned}
\end{equation}
where \(s_{\lambda}\) represents a Schur function associated to a partition \(\lambda\).
\end{lem}

It can be easily deduced from the following proposition from Stanley [12].

\begin{prop}
\([\)12, (7.105)\(]\) For any set \(X=\{x_1,\dotsc,x_n\}\) of positive integers, where elements are written in increasing order, we have
\begin{equation}
    s_{\lambda(X)}(1,q,\dotsc,q^{n-1})=\frac{\Delta_{1,q}(X)}{\Delta_{1,q}([n])}.
\end{equation}
\end{prop}

\textit{Proof of Lemma 3.2.}
By Proposition 3.3, for any such set \(Z\), we have
\begin{equation}
\begin{aligned}
    &\frac{s_{\lambda(X'\cup Z)}(1,q,\dotsc,q^{x'+z-1})s_{\lambda(Y'\cup Z)}(1,q,\dotsc,q^{y'+z-1})}{s_{\lambda(X\cup Z)}(1,q,\dotsc,q^{x+z-1})s_{\lambda(Y\cup Z)}(1,q,\dotsc,q^{y+z-1})}\\
    &=\frac{\Delta_{1,q}([x+z])\Delta_{1,q}([y+z])}{\Delta_{1,q}([x'+z])\Delta_{1,q}([y'+z])}\frac{\Delta_{1,q}(X'\cup Z)\Delta_{1,q}(Y'\cup Z)}{\Delta_{1,q}(X\cup Z)\Delta_{1,q}(Y\cup Z)}.
\end{aligned}
\end{equation}
We can simplify the terms containing \(Z\) on the right side of (3.4) as follows:

\begin{equation}
\begin{aligned}
    \frac{\Delta_{1,q}(X'\cup Z)\Delta_{1,q}(Y'\cup Z)}{\Delta_{1,q}(X\cup Z)\Delta_{1,q}(Y\cup Z)}
    &=\frac{\Delta_{1,q}(X')\Delta_{2,q}(X', Z)\Delta_{1,q}(Z)\Delta_{1,q}(Y')\Delta_{2,q}(Y', Z)\Delta_{1,q}(Z)}{\Delta_{1,q}(X)\Delta_{2,q}(X, Z)\Delta_{1,q}(Z)\Delta_{1,q}(Y)\Delta_{2,q}(Y, Z)\Delta_{1,q}(Z)}\\
    &=\frac{\Delta_{1,q}(X')\Delta_{1,q}(Y')}{\Delta_{1,q}(X)\Delta_{1,q}(Y)}.
\end{aligned}
\end{equation}

Hence, by (3.4) and (3.5),
\begin{equation}
\begin{aligned}
    &\frac{s_{\lambda(X'\cup Z)}(1,q,\dotsc,q^{x'+z-1})s_{\lambda(Y'\cup Z)}(1,q,\dotsc,q^{y'+z-1})}{s_{\lambda(X\cup Z)}(1,q,\dotsc,q^{x+z-1})s_{\lambda(Y\cup Z)}(1,q,\dotsc,q^{y+z-1})}\\
    &=\frac{\Delta_{1,q}([x+z])\Delta_{1,q}([y+z])}{\Delta_{1,q}([x'+z])\Delta_{1,q}([y'+z])}\frac{\Delta_{1,q}(X')\Delta_{1,q}(Y')}{\Delta_{1,q}(X)\Delta_{1,q}(Y)}.
\end{aligned}
\end{equation}
This completes the proof. \(\square\)

\medskip
The identity above will be used in the proof of Theorem 2.1 (a) via the following well-known relation between the Schur function and weighted enumeration of lozenge tilings of the trapezoidal region with some dents. For reference, we state the version of Ayyer and Fischer from [1] (in [1], they stated it in terms of matching generating function of a certain graph, which is equivalent).

\begin{thm}
\([\)1, Theorem 2.3\(]\) Consider the region \(T_{m,n}(S)\) that we described in Proposition 3.1. On this region, we give weight \(t_k\) to each right-lozenge whose distance between the bottom side of the lozenge and the top side of the region is \(\frac{k\sqrt{3}}{2}\), and give weight 1 to all vertical- and left-lozenges. If \(M(T_{m,n}(S);(t_1,\dotsc,t_n))\) is the tiling generating function of \(T_{m,n}(S)\) under this weight, then we have
\begin{equation}
    M(T_{m,n}(S);(t_1,\dotsc,t_n))=s_{\lambda(S)}(t_1,\dotsc,t_n).
\end{equation}
\end{thm}

Recall that Schur functions are symmetric and homogeneous. Theorem 3.4 allows us to convert these properties of Schur functions into the following properties of tiling generating functions:
\begin{equation*}
\begin{aligned}
    &M(T_{m,n}(S);(t_1,\dotsc,t_n))=M(T_{m,n}(S);(t_{\sigma(1)},\dotsc,t_{\sigma(n)})), \forall \sigma\in S_n \text{, and}\\
    &M(T_{m,n}(S);(qt_1,\dotsc,qt_n))=q^{|\lambda(S)|}M(T_{m,n}(S);(t_1,\dotsc,t_n)).
\end{aligned}
\end{equation*}

Note that for any bounded region \(G\) whose height is \(\frac{n\sqrt{3}}{2}\), we have 
\begin{equation*}
    M(G;(q,q^2,\dotsc,q^n))=M(G;q).
\end{equation*}

By using Lemma 3.2 and the properties of tiling generating functions, we can give a simple proof of Theorem 2.1 (a).

To prove part (b) (the shuffling theorem for centrally symmetric tilings), we need the following simple lemma.

\begin{lem}
Let \(k\) be a positive integer, and let \(X\), \(X'\)  be subsets of \([k]\) such that a pair \((X', X'_{(k)})\) is a shuffling of \((X, X_{(k)})\). If \(Z\) is any subset of \([k]\) with cardinality \(z\) that is disjoint from \(X\cup X_{(k)} (=X'\cup X'_{(k)})\) and is \(k\)-symmetric with itself \((Z=Z_{(k)})\), then we have
\begin{equation}
    \frac{\Delta_{1}(X'\cup Z)}{\Delta_{1}(X\cup Z)}=\frac{\Delta_{1}(X')}{\Delta_{1}(X)}.
\end{equation}
\end{lem}

The proof of Lemma 3.5 is analogous to that of Lemma 3.2. Additionally, we have to use the \(k\)-symmetric relations between some sets.

\medskip
\textit{Proof of Lemma 3.5.}
By shuffling condition and \(k\)-symmetric relations, we have
\begin{equation}
\begin{aligned}
    \Delta_{2}(X,Z)=\sqrt{\Delta_{2}(X,Z)\cdot\Delta_{2}(X_{(k)},Z_{(k)})}&=\sqrt{\Delta_{2}(X,Z)\cdot\Delta_{2}(X_{(k)},Z)}\\
    &=\sqrt{\Delta_{2}(X',Z)\cdot\Delta_{2}(X'_{(k)},Z)}\\
    &=\sqrt{\Delta_{2}(X',Z)\cdot\Delta_{2}(X'_{(k)},Z_{(k)})}\\
    &=\Delta_{2}(X',Z).\\
\end{aligned}
\end{equation}

By (3.9), we have
\begin{equation}
\begin{aligned}
    \frac{\Delta_{1}(X'\cup Z)}{\Delta_{1}(X\cup Z)}=\frac{\Delta_{1}(X') \Delta_{2}(X',Z) \Delta_{1}(Z)}{\Delta_{1}(X) \Delta_{2}(X,Z) \Delta_{1}(Z)}=\frac{\Delta_{1}(X')}{\Delta_{1}(X)}.
\end{aligned}
\end{equation}
This completes the proof. \(\square\)

Using Proposition 3.1 and Lemma 3.5, we will provide a short proof of Theorem 2.1 (b).

\section{The proof of Theorem 2.1}

We first prove part (a). The set of lozenge tilings of \(V_{a, b, c}(X,Y)\) can be partitioned according to the positions of \(b-x\) vertical-lozenges crossing the horizontal diagonal \(\ell\). Let \(Z\) be an index set for positions of the vertical-lozenges. If we fix these vertical-lozenges, then the remaining region to be tiled is \(V_{a, b, c}(X\cup Z,Y\cup Z)\). Furthermore, since tilings of \(V_{a, b, c}(X\cup Z,Y\cup Z)\) have no vertical-lozenges along \(\ell\), tilings of \(V_{a, b, c}(X\cup Z,Y\cup Z)\) are in bijection with pairs of tilings, one from tilings of \(T_{a,b}(X\cup Z)\) and the other from tilings of \(T_{a+b-c,c}(Y\cup Z)_{(a+b)}\), where \((Y\cup Z)_{(a+b)}=\{a+b+1-y|y\in Y\cup Z\}\) (see Figure 4.1). Hence,

\begin{equation}
    M(V_{a,b,c}(X,Y);q)=\sum_{Z}M(V_{a,b,c}(X\cup Z,Y\cup Z);q),
\end{equation}
where \(Z\) runs over all subsets of \([a+b]\setminus(X\cup Y)\) with cardinality \(b-x\),
and
\begin{equation}
    M(V_{a,b,c}(X\cup Z,Y\cup Z);q)=M(T_{a,b}(X\cup Z);q) M(T_{a+b-c,c}(Y\cup Z)_{(a+b)};(q^{b+c},q^{b+c-1},\dotsc,q^{b+1})).
\end{equation}

\begin{figure}
    \centering
    \includegraphics[width=10cm]{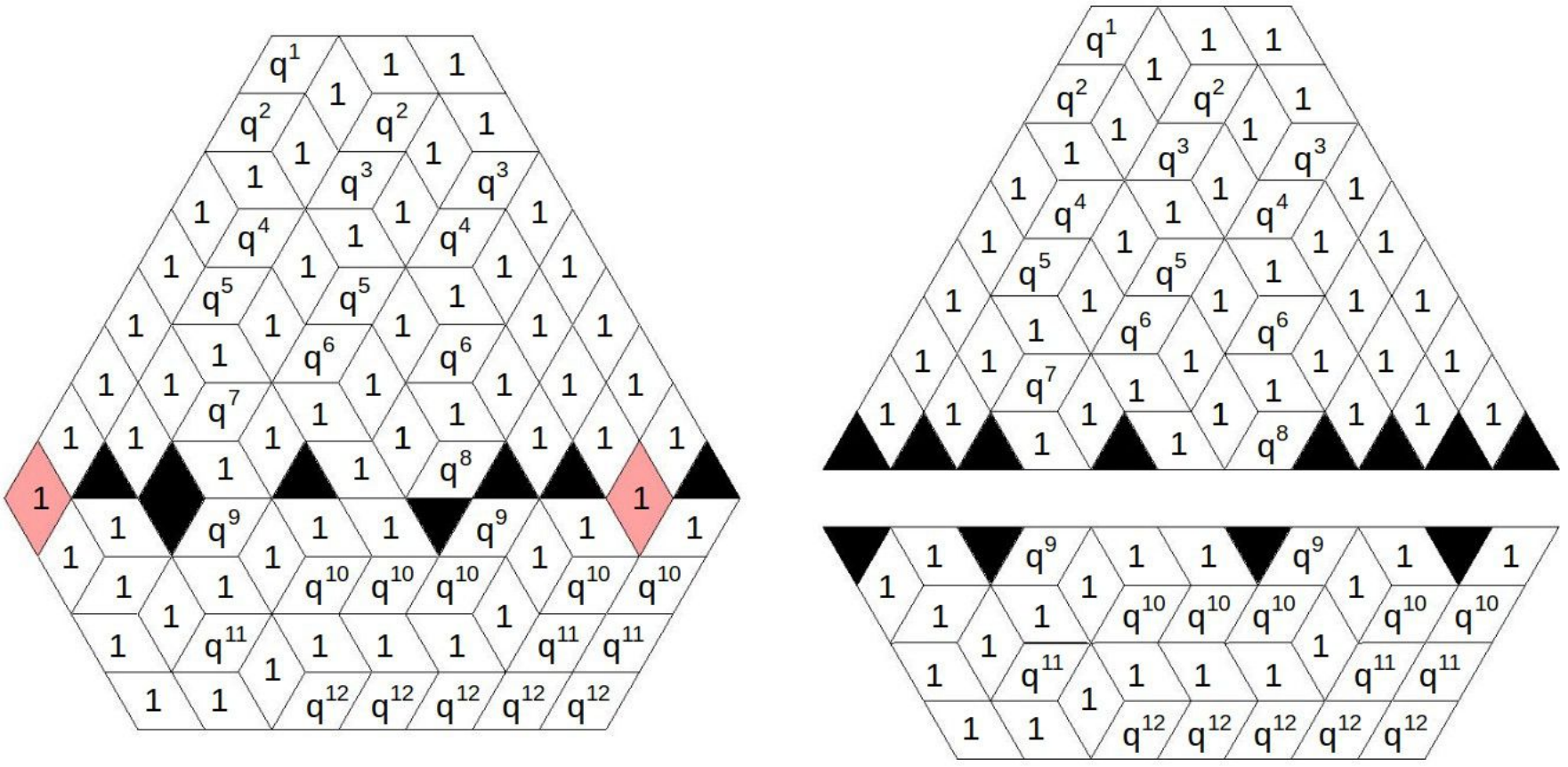}
    \caption{A lozenge tiling of \(V_{3,8,4}(\{2,3,5,8,9,11\},\{3,7\})\) with vertical--lozenges at positions \(\{1,10\}\) (left) and corresponding pair of lozenge tilings of two subregions (right). Note that each lozenge has the weight.}
\end{figure}

Also, by Theorem 3.4 and the properties of tiling generating functions, 
\begin{equation}
\begin{aligned}
    M(T_{a,b}(X\cup Z);q)&=M(T_{a,b}(X\cup Z);(q,q^2,\dotsc,q^{b}))\\
    &=q^{|\lambda(X\cup Z)|}M(T_{a,b}(X\cup Z);(1,q,\dotsc,q^{b-1}))\\
    &=q^{(\sum_{i\in X}i+\sum_{k\in Z}k-\sum_{l=1}^{b}l)}s_{\lambda(X\cup Z)}(1,q,\dotsc,q^{b-1})
\end{aligned}
\end{equation}
and similarly
\begin{equation}
\begin{aligned}
    &M(T_{a+b-c,c}(Y\cup Z)_{(a+b)};(q^{b+c},q^{b+c-1},\dotsc,q^{b+1}))\\
    &=M(T_{a+b-c,c}(Y\cup Z)_{(a+b)};(q^{b+1},q^{b+2},\dotsc q^{b+c}))\\
    &=q^{(b+1)|\lambda((Y\cup Z)_{(a+b)})|}M(T_{a+b-c,c}(Y\cup Z)_{(a+b)};(1,q,\dotsc q^{c-1}))\\
    &=q^{(b+1)\{(a+b+1)c-\sum_{j\in Y}j-\sum_{k\in Z}k-\sum_{l=1}^{c}l\}}s_{\lambda((Y\cup Z)_{(a+b)})}(1,q,\dotsc,q^{c-1}).
\end{aligned}
\end{equation}

By Proposition 3.3, we have
\begin{equation}
\begin{aligned}
    s_{\lambda((Y\cup Z)_{(a+b)})}(1,q,\dotsc,q^{c-1})
    &=\frac{\Delta_{1,q}((Y\cup Z)_{(a+b)})}{\Delta_{1,q}([c])}\\
    &=\frac{1}{\Delta_{1,q}([c])}\displaystyle \prod_{i, j \in (Y\cup Z)_{(a+b)}, i<j}([j]_q-[i]_q)\\
    &=\frac{1}{\Delta_{1,q}([c])}\displaystyle \prod_{i, j \in Y\cup Z, j<i}([a+b+1-j]_q-[a+b+1-i]_q)\\
    &=\frac{1}{\Delta_{1,q}([c])}\displaystyle \prod_{i, j \in Y\cup Z, j<i}q^{a+b+1-i-j}([i]_q-[j]_q)\\
    &=\Bigg(\prod_{i, j \in Y\cup Z, j<i}q^{a+b+1-i-j}\Bigg)  \frac{\Delta_{1,q}(Y\cup Z)}{\Delta_{1,q}([c])}\\
    &=q^{\{(a+b+1)\frac{c(c-1)}{2}-(c-1)(\sum_{j\in Y}j+\sum_{k\in Z}k)\}}s_{\lambda(Y\cup Z)}(1,q,\dotsc,q^{c-1}).\\
\end{aligned}
\end{equation}
Thus, by (4.2)-(4.5),
\begin{equation}
    M(V_{a,b,c}(X\cup Z,Y\cup Z);q)=q^{\alpha(Z)} s_{\lambda(X\cup Z)}(1,q,\dotsc,q^{b-1})s_{\lambda(Y\cup Z)}(1,q,\dotsc,q^{c-1}),
\end{equation}
where \(\displaystyle\alpha(Z)=\sum_{i\in X}i-(b+c)\sum_{j\in Y}j+(1-b-c)\sum_{k\in Z}k-\frac{b(b+1)}{2}+(a+b+1)(b+1)c-\frac{(b+1)c(c+1)}{2}+\frac{(a+b+1)c(c-1)}{2}\).

Similarly, we also have
\begin{equation}
    M(V_{a',b' ,c'}(X',Y');q)=\sum_{Z}M(V_{a',b',c'}(X'\cup Z,Y'\cup Z);q),
\end{equation}
where \(Z\) runs over all subsets of \([a'+b']\setminus(X'\cup Y') (= [a+b]\setminus(X\cup Y))\) whose cardinality is \(b'-x' (= b-x)\), and
\begin{equation}
    M(V_{a',b',c'}(X'\cup Z,Y'\cup Z);q)=q^{\alpha'(Z)} s_{\lambda(X'\cup Z)}(1,q,\dotsc,q^{b-1})s_{\lambda(Y'\cup Z)}(1,q,\dotsc,q^{c-1}),
\end{equation}
where \(\displaystyle\alpha'(Z)=\sum_{i'\in X'}i'-(b'+c')\sum_{j'\in Y'}j'+(1-b'-c')\sum_{k\in Z}k-\frac{b'(b'+1)}{2}+(a'+b'+1)(b'+1)c'-\frac{(b'+1)c'(c'+1)}{2}+\frac{(a'+b'+1)c'(c'-1)}{2}\).

One can readily check that summations in (4.1) and (4.7) are taken over the same sets \(Z\).
For any such \(Z\), by (4.6), (4.8) and Lemma 3.2, the ratio of the corresponding summands is

\begin{equation}
\begin{aligned}
    \frac{M(V_{a',b',c'}(X'\cup Z,Y'\cup Z);q)}{M(V_{a,b,c}(X\cup Z,Y\cup Z);q)}
    &=q^{\alpha'(Z)-\alpha(Z)}\frac{s_{\lambda(X'\cup Z)}(1,q,\dotsc,q^{b'-1})s_{\lambda(Y'\cup Z)}(1,q,\dotsc,q^{c'-1})}{s_{\lambda(X\cup Z)}(1,q,\dotsc,q^{b-1})s_{\lambda(Y\cup Z)}(1,q,\dotsc,q^{c-1})}\\
    &=q^{\alpha}\frac{\Delta_{1,q}([b])\Delta_{1,q}([c])}{\Delta_{1,q}([b'])\Delta_{1,q}([c'])}\frac{\Delta_{1,q}(X')\Delta_{1,q}(Y')}{\Delta_{1,q}(X)\Delta_{1,q}(Y)},
\end{aligned}
\end{equation}
where \(\displaystyle\alpha=\alpha'(Z)-\alpha(Z)=\Bigg[\sum_{i'\in X'}i'-(b'+c')\sum_{j'\in Y'}j'-\frac{b'(b'+1)}{2}+(a'+b'+1)(b'+1)c'-\frac{(b'+1)c'(c'+1)}{2}+\frac{(a'+b'+1)c'(c'-1)}{2}\Bigg]-\Bigg[\sum_{i\in X}i-(b+c)\sum_{j\in Y}j-\frac{b(b+1)}{2}+(a+b+1)(b+1)c-\frac{(b+1)c(c+1)}{2}+\frac{(a+b+1)c(c-1)}{2}\Bigg]\).

The expression on the right side of (4.9) does not depend on the set \(Z\).
Therefore, by (4.1), (4.7), and (4.9),
\begin{equation}
    \frac{M(V_{a',b',c'}(X',Y');q)}{M(V_{a,b,c}(X,Y);q)}=q^{\alpha}\frac{\Delta_{1,q}([b])\Delta_{1,q}([c])}{\Delta_{1,q}([b'])\Delta_{1,q}([c'])}\frac{\Delta_{1,q}(X')\Delta_{1,q}(Y')}{\Delta_{1,q}(X)\Delta_{1,q}(Y)}.
\end{equation}
This completes the proof of part (a).

Using the same idea, now we prove part (b). Again, we partition the set of centrally symmetric lozenge tilings of the region \(V_{a,b,b}(X,X_{(a+b)})\) according to \((b-x)\) vertical-lozenges crossing the horizontal diagonal. If \(Z\) is the index set that record the positions of the vertical-lozenges, then this set \(Z\) should be \((a+b)\)-symmetric with itself (\(Z=Z_{(a+b)}\)). If we fix these vertical-lozenges, then the remaining region to be tiled is \(V_{a, b, b}(X\cup Z,X_{(a+b)}\cup Z_{(a+b)})\) and one can readily see that its centrally symmetric lozenge tilings are in bijection with lozenge tilings of its subregion above the horizontal diagonal, \(T_{a,b}(X\cup Z)\) (see Figure 4.2). Hence,
\begin{equation}
    M_\odot(V_{a, b, b}(X,X_{(a+b)}))=\sum_{Z}M(T_{a,b}(X\cup Z))=\frac{\displaystyle\sum_{Z}\Delta_{1}(X\cup Z)}{H(b)},
\end{equation}
where the sum is taken over all \((b-x)\)-element sets \(Z\subset[a+b]\setminus(X\cup X_{(a+b)})\) such that \(Z=Z_{(a+b)}\) holds.
By exactly the same argument,
\begin{equation}
    M_\odot(V_{a, b, b}(X',X'_{(a+b)}))=\sum_{Z}M(T_{a,b}(X'\cup Z))=\frac{\displaystyle\sum_{Z}\Delta_{1}(X'\cup Z)}{H(b)},
\end{equation}
where the sum is taken over all \((b-x)\)-element sets \(Z\subset[a+b]\setminus(X'\cup X'_{(a+b)})(=[a+b]\setminus(X\cup X_{(a+b)}))\)  such that \(Z=Z_{(a+b)}\) holds.

\medskip
\begin{figure}
    \centering
    \includegraphics[width=10cm]{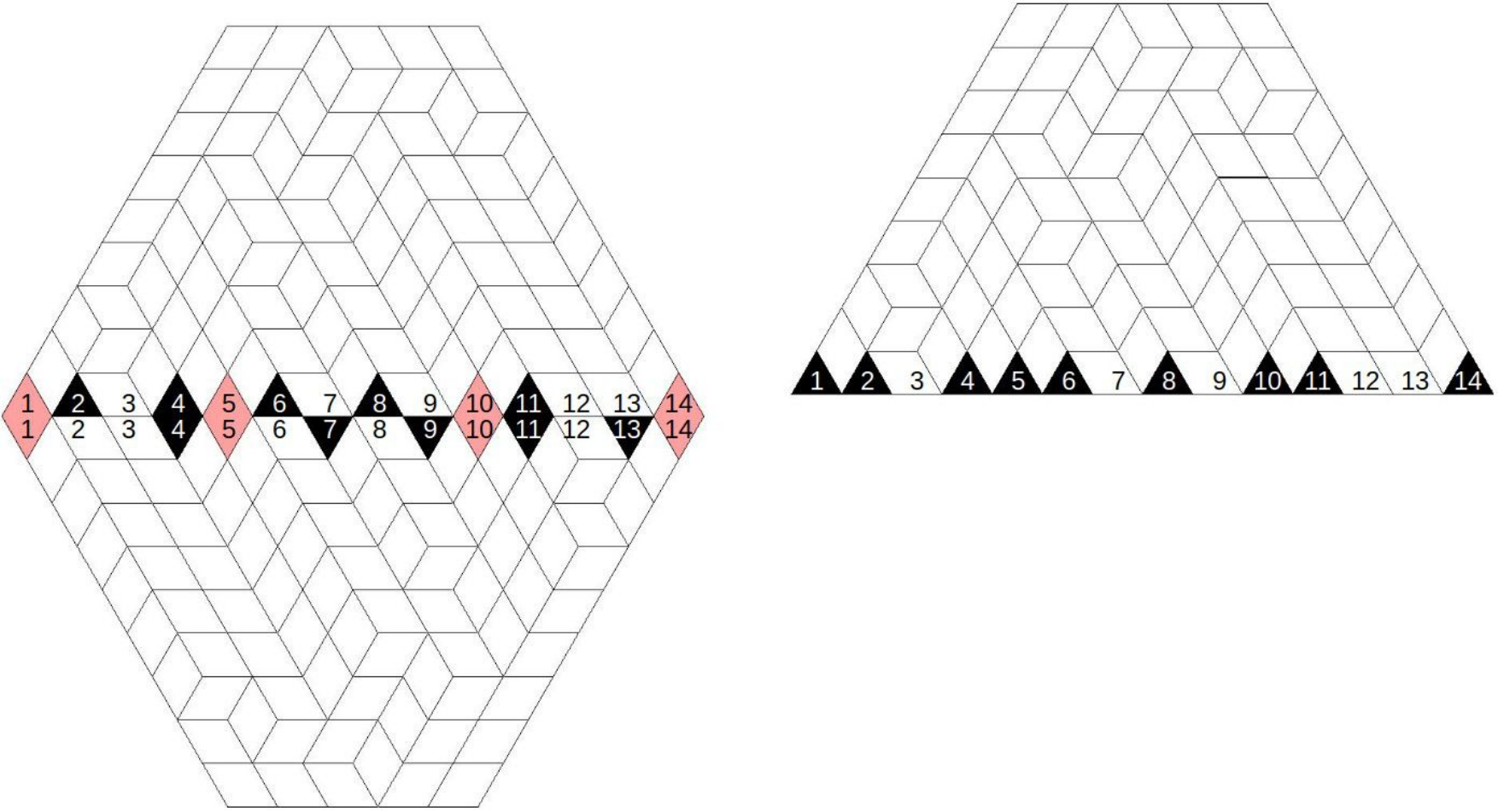}
    \caption{A centrally symmetric lozenge tiling of \(V_{5,9,9}(\{2,4,6,8,11\},\{4,7,9,11,13\})\) with vertical-lozenges at positions \(\{1, 5, 10, 14\}\) (left) and the corresponding lozenge tiling of its subregion \(T_{5,9}(\{1,2,4,5,6,8,10,11,14\})\) (right).}
\end{figure}

Note that summations in (4.11) and (4.12) are taken over the same sets \(Z\).
For any such \(Z\), Lemma 3.5 says the ratio of the corresponding summands is
\begin{equation}
    \frac{\Delta_{1}(X'\cup Z)}{\Delta_{1}(X\cup Z)}=\frac{\Delta_{1}(X')}{\Delta_{1}(X)}.
\end{equation}

This ratio does not depend on the set \(Z\). Thus, by (4.11)-(4.13), we have 
\begin{equation}
    \frac{M_\odot(V_{a, b, b}(X',X'_{(a+b)}))}{M_\odot(V_{a, b, b}(X,X_{(a+b)}))}= \frac{\displaystyle\sum_{Z}M(T_{a,b}(X'\cup Z))}{\displaystyle\sum_{Z}M(T_{a,b}(X\cup Z))}=\frac{\Delta_{1}(X')}{\Delta_{1}(X)}.
\end{equation}
Remaining part of the theorem is clear from Theorem 2.1 (a) (with $q=1$) and the two facts \(\Delta_{1}(X)=\Delta_{1}(X_{(a+b)})\) and \(\Delta_{1}(X')=\Delta_{1}(X'_{(a+b)})\), which can be easily deduced from \((a+b)\)-symmetric relations between sets.
This completes the proof. \(\square\)

\medskip
It is immediate from the proof that the ratio between the number of centrally symmetric tilings of the two regions is equal to the square root of that of the total number of tilings. This is because centrally symmetric tilings are determined by the lozenges above the horizontal line, while both lozenges above and below the horizontal line contribute in the other case.

\section{Acknowledgment}
The author would like to thank his advisor, Professor Mihai Ciucu, for his encouragement and useful discussions. This paper could not be written without his continued guidance. The author also thanks anonymous referees for carefully reading the earlier version of the paper and giving helpful feedback. Also, the author thanks Jeff Taylor for installing software and helpful assistance. David Wilson's program, \textit{vaxmacs}, was extremely useful when the author made an observation.

\end{document}